\documentclass[a4paper,11pt]{article}
\usepackage[usenames]{color}
\usepackage[utf8]{inputenc}
\usepackage{amsmath}
\usepackage{amssymb}
\usepackage{amsfonts}
\usepackage{amsthm,amscd}
\usepackage[utf8]{inputenc}
\newtheorem{theorem}{Theorem}[section]
\newtheorem{lemma}{Lemma}[section]
\newtheorem{corollary}{Corollary}[section]
\newtheorem{proposition}{Proposition}[section]
\theoremstyle{definition}
\newtheorem{definition}{Definition}[section]

\theoremstyle{remark}
\newtheorem{remark}[theorem]{Remark}

\numberwithin{equation}{section}

\title{On the geometry of co-Hamiltonian  diffeomorphisms}
\author{\scshape 
	S. Tchuiaga	$^{}$\thanks{tchuiagas@gmail.com,\ Department of Mathematics,  
		The University of Buea, South West Region,Cameroon},  and P. Bikorimana	$^{}$\thanks{pierrebikorimana@gmail.com, \ Institut de Mathematiques et des Science Physiques, Benin} }

\definecolor{couleurliens}{rgb}{1.0,0.,0.} 
 \definecolor{couleurliensref}{rgb}{0.,0.,1.} 
\definecolor{couleurliensurl}{rgb}{.3,.4,.3} 
\usepackage[pdftex,colorlinks=true,urlcolor=couleurliensurl,linkcolor=couleurliens,citecolor=couleurliensref]{hyperref}

\begin{document}
\maketitle\large 
\begin{abstract}  
This paper studies the geometry of the group of all co-Hamiltonian diffeomorphisms  of a compact  cosymplectic manifold $(M, \omega, \eta)$. The 
fix-point theory for co-Hamiltonian diffeomorphisms is studied, and we use Arnold's conjecture to predict the exact 
minimum number of fix point that such a diffeomorphism must have (this minimum number is at least $1$).  It follows 
that the generating function of any co-Hamiltonian isotopy is a constant function along it orbits. 
Therefore, we study the co-Hofer norms for co-Hamiltonian isotopies, and 
establish several co-Hamiltonian and almost co-Hamiltonian analogues of some approximations lemmas and reparameterizations lemmas
found in the theory of Hamiltonian dynamics, we define two $C^0-$co-Hamiltonian topologies, 
and use these topologies to define the spaces of cohameomorphisms, and almost cohameomorphisms. 
Finally,  we raise several important questions for future studies. 
 
\end{abstract}


\vspace{1cm}

\section{Preliminaries}\label{SC0}
\subsection{Cosymplectic vector spaces}
A bilinear form on a vector space $V$ is a map $b : V\times V \longrightarrow \mathbb{R}$ which is linear in each variable. When $b(u, v) = - b(v,u)$, for all $u, v\in V$, then $b$ is called antisymmetric.\\

Given any non-trivial linear map $L:V\longrightarrow \mathbb{R}$, together with a bilinear map\\ $b : V\times V \longrightarrow \mathbb{R}$, one defines a linear map $$\tilde{I}_{L, b}: V\longrightarrow V^\ast, Y\mapsto \iota(Y)b + L(Y)L,$$
so that $\tilde{I}_{L, b}(Y)(X)  = b(Y,X) + L(Y)L(X)$, for all $X, Y \in V$, where $V^\ast$ is the dual space of $V$.\\

\begin{definition}\cite{T-H-B}
A pair $(b, L)$ consists of  antisymmetric bilinear map $b : V\times V \longrightarrow \mathbb{R}$, and non-trivial linear map $L:V\longrightarrow \mathbb{R}$ is called 
 cosymplectic couple if the map  $ \tilde{I}_{L, b}$ is a bijection. 
\end{definition}
\begin{definition}\cite{T-H-B}
	A cosymplectic vector space is a triple $(V, b, L)$ where  $ V$ is a vector space and $(b, L)$ is a  
	cosymplectic couple. 
\end{definition}
\begin{proposition}\label{Dim}
	Let $(V, b, L)$ be a cosymplectic vector space. Then, $\dim (V) = 2n + 1$. 
\end{proposition}

\begin{remark}\label{Dim1}	Let $(V, b, L)$ be a cosymplectic vector space. From the proof of Proposition \ref{Dim}, one can always assume that $ L(\xi) = 1$ (of course after normalization, if necessary), and we also have $b(\xi, v) = 0$, for all $v\in V$. We shall call the vector $\xi$, the Reeb vector of  $(V, b, L)$. 
\end{remark}

\subsection{Cosymplectic Manifolds, \cite{T-H-B}}
Let $M$ be a smooth  manifold. An almost cosymplectic structure on $M$ is a pair $(\eta, \omega)$ consisting
of a  $1-$form $\eta$, and $2-$form $\omega$ such that for each $x\in M$, the triple $ (T_xM, \omega_x, \eta_x)$ is a cosymplectic vector space.  Therefore, a cosymplectic structure on $M$ is any 
almost cosymplectic structure $(\eta, \omega)$ on $M$ such that $d\eta = 0$, and $d\omega = 0$. We shall write  $(M, \eta, \omega)$ to mean that $M$ is equipped with a cosymplectic structure $(\eta, \omega)$. In particular, any cosymplectic manifold $(M, \eta, \omega)$ is oriented with respect to the volume form 
$\eta\wedge\omega^n$, while any cosymplectic manifold $(M, \eta, \omega)$ admits a vector field $\xi$ called the 
Reeb vector field such that $ \eta(\xi) = 1$, and $\iota(\xi)\omega = 0$. We refer the readers to the papers \cite{P.L, H-L, T-H-B} for further details.\\ 

\begin{definition}\cite{T-H-B}
	 Let $(M, \eta, \omega)$ be a cosymplectic manifold. A diffeomorphism $\phi : M\longrightarrow M$ is said to be a  cosymplectic diffeomorphism (or cosymplectomorphism) if: $\phi^\ast(\eta)= \eta$, and $\phi^\ast(\omega)= \omega$. 
\end{definition}
We shall denote by $Cosymp_{\eta, \omega}^\ast(M)$ the space of all cosymplectomorphisms of $(M, \eta, \omega)$. 
\begin{definition}\cite{T-H-B} Let $(M, \eta, \omega)$ be a cosymplectic manifold. An isotopy $\Phi =\{\phi_t\}$ is said to be cosymplectic if for each $t$, we have 
 $\phi_t\in  Cosymp_{\eta, \omega}^\ast(M)$. 
\end{definition}
This group is shall denote by  $Iso_{\eta, \omega}^\ast(M)$ the space of all cosymplectic isotopies   of $(M, \eta, \omega)$, and put
\begin{equation}
G_{\eta, \omega}^\ast(M) := ev_1\left(Iso_{\eta, \omega}^\ast(M)\right). 
\end{equation}
We equip the group $G_{\eta, \omega}^\ast(M) $ with the $C^\infty-$compact-open topology \cite{Hirs76}.
\begin{definition}\cite{T-H-B}
	Let $(M, \eta, \omega)$ be a cosymplectic manifold. An element $Y\in \chi_{\eta, \omega}(M)$, is called an almost  co-Hamiltronian vector field if 
	the $1-$form $\iota(\omega)(Y)$ is exact. 
\end{definition}
We shall denote by $ham_{\eta, \omega}(M)$ the space of all almost co-Hamiltonian vector fields of $(M, \eta, \omega)$. 
\begin{definition}\cite{T-H-B}
	Let $(M, \eta, \omega)$ be a compact cosymplectic manifold. A cosymplectic isotopy $\Psi:=\{\psi_t\}$ is called an almost co-Hamiltonian isotopy, if
	for each $t$, the vector field $\dot\psi_t$ is a co-Hamiltronian vector field, i.e., $\dot\psi_t\in ham_{\eta, \omega}(M)$, for each $t$. 
\end{definition}
We shall denote by $\mathcal{A}H_{\eta, \omega}(M)$ the space of all almost co-Hamiltonian isotopies of $(M, \eta, \omega)$, and put
\begin{equation}
\mathcal{A}Ham_{\eta, \omega}(M) := ev_1\left(\mathcal{A}H_{\eta, \omega}(M)\right). 
\end{equation}
The elements of the set $\mathcal{A}H_{\eta, \omega}(M)$ are called almost co-Hamiltonian diffeomorphisms of $(M, \eta, \omega)$.

\subsubsection{Cosymplectic vector fields, \cite{T-H-B}}\label{SC011}

\begin{definition}
	Let $(M, \eta, \omega)$ be a  cosymplectic manifold. A vector field $X$ is said to be  cosymplectic if $\mathcal{L}_{X}\eta = 0$,
	and  $\mathcal{L}_{X}\omega = 0$.
\end{definition}
We shall denote by $\chi_{\eta, \omega}^\ast(M)$ the space of all  cosymplectic vector fields of $(M, \eta, \omega)$.
\begin{definition}
	Let $(M, \eta, \omega)$ be a cosymplectic manifold. An element $Y\in \chi_{\eta, \omega}^\ast(M)$, is called a co-Hamiltonian vector field if 
	the $1-$form $I_{\eta, \omega}(Y)$ is exact. 
\end{definition}
We shall denote by $ham_{\eta, \omega}^\ast(M)$ the space of all co-Hamiltonian vector fields of $(M, \eta, \omega)$. 
 
\begin{definition}
	Let $(M, \eta, \omega)$ be a compact cosymplectic manifold. A cosymplectic isotopy $\Psi:=\{\psi_t\}$ is called a co-Hamiltonian isotopy, if
	for each $t$, the vector field $\dot\psi_t$ is a co-Hamiltonian vector field, i.e., $\dot\psi_t\in ham_{\eta, \omega}^\ast(M)$, for each $t$. 
\end{definition}
We shall denote by $H_{\eta, \omega}(M)$ the space of all  co-Hamiltonian isotopies of $(M, \eta, \omega)$, and put
\begin{equation}
Ham_{\eta, \omega}(M) := ev_1\left(H_{\eta, \omega}(M)\right). 
\end{equation}
The elements of the set $Ham_{\eta, \omega}(M)$ are called co-Hamiltonian diffeomorphisms of $(M, \eta, \omega)$. 
\begin{proposition}\cite{T-H-B}
	Let $(M, \eta, \omega)$ be a cosymplectic manifold. For any $X, Y\in \chi_{\eta, \omega}^\ast(M)$, we have $[X,Y]\in ham_{\eta, \omega}^\ast(M)$.
\end{proposition}
We shall need the following facts that we borrow from \cite{T-H-B}. 
	\begin{enumerate}
		\item Let $\{\phi_t\}\in H_{\eta, \omega}(M)$ such that $I_{\eta, \omega}(\dot\phi_t) = df_t$, for each $t$. Then, $\{\phi_t^{-1}\}\in H_{\eta, \omega}(M)$, and $I_{\eta, \omega}(\dot{\phi}_{-t}) =  d(-f_t\circ\phi_t), $ where 
		$\phi_t^{-1} =:\phi_{-t}$, for each $t$. 
		\item If  $\Phi_F = \{\phi_t\}$ is a co-Hamiltonian isotopy such that $I_{\eta, \omega} (\dot\phi_t) = dF_t$, for all $t$, then for all 
		$\rho \in G_{\eta,  \omega}^\ast(M)$, the isotopy $\bar \Phi : t\mapsto \rho^{-1} \circ \phi_t\circ \rho$ is co-Hamiltonian, and it is generated by 
		$ H_t :=  F_t\circ\rho,$ 
		for each $t$. 
		\item If  $\{\psi_t\}, \{\phi_t\}\in H_{\eta, \omega}(M)$ such that $I_{\eta, \omega}(\dot\phi_t) = df_t$, and 
		$I_{\eta, \omega}(\dot\psi_t) = dh_t$,
		for each $t$, then the isotopy $t\mapsto \phi_t\circ \phi_t$ is generated by, $h_t = f_t + h_t\circ\phi_t^{-1}$, for each $t$.
		\item Let $\{\phi_t\}\in Iso_{\eta, \omega}^\ast(M)$. Then,
		for each $t$, we have  
		$I_{\eta, \omega}(\dot{\phi}_{-t}) = -\phi_t^\ast(I_{\eta, \omega}(\dot{\phi}_t)) .$
		\item If  $\{\psi_t\}, \{\phi_t\}\in Iso_{\eta, \omega}^\ast(M)$, then 
		for each $t$, we have 
		$$I_{\eta, \omega}(\dot{\overline{\phi_t\circ\psi_t}}) = I_{\eta, \omega}(\dot{\phi}_t) + (\phi_t^{-1})^\ast(I_{\eta, \omega}(\dot{\psi}_t)).$$
		\item Let $\{\phi_t\}\in Iso_{\eta, \omega}(M) $ such that $\mathcal{L}_{\dot\phi_t}\eta = \mu_t\eta$, (or $\phi_t^\ast(\eta) = e^{f_t}\eta$), 
		for each $t$. We have  $\mathcal{L}_{\dot{\phi}_{-t}}\eta = \vartheta_t\eta,$ for all $t$, with 
		$ \vartheta_t := -(\dot{\overline{f_t\circ\phi_t^{-1}}})\circ\phi_t.$ 
		\item If  $\{\psi_t\}, \{\phi_t\}\in Iso_{\eta, \omega}(M)$
		such that $\mathcal{L}_{\dot\phi_t}\eta = \mu_t\eta$, (or $\phi_t^\ast(\eta) = e^{f_t}\eta$), and 
		$\mathcal{L}_{\dot\psi_t}\eta = \mu'_t\eta$, (or $\psi_t^\ast(\eta) = e^{q_t}\eta$),
		for each $t$, then we have $ 	\mathcal{L}_{\dot{\overline{\phi_t\circ\psi_t}}}\eta = \varrho_t\eta,$ 
		with\\  $ \varrho_t : = \left( \dot{\overline{f_t\circ\psi_t}} + \dot{q}_t \right)\circ(\phi_t\circ\psi_t)^{-1},$ 
		for all $t$.
		\item 	If  $\Phi = \{\phi_t\}$ is an almost co-Hamiltonian isotopy such that $ \iota(\dot\phi_t)\omega = dF_t$, for all $t$, and \\
		$\mathcal{L}_{\dot\phi_t}\eta = \mu_t\eta$, (or $\phi_t^\ast(\eta) = e^{f_t}\eta$), 
		for each $t$, then  for all $\rho \in G_{\eta,  \omega}(M)$ such that\\ $\rho^\ast(\eta) = e^{f^\rho}\eta$, 
		the isotopy $\bar \Phi : t\mapsto \rho^{-1} \circ \phi_t\circ \rho$ is an almost co-Hamiltonian, and we have 
		$ \mathcal{L}_{\rho_\ast^{-1}(\dot\phi_t)}\eta =  H_t(\rho)\eta,$
		with 
		$ H_t(\rho)  = \left(-df^\rho((\rho^{-1})_\ast(\dot{\phi}_t )) + \dot f_t\right) \circ (\rho^{-1}\circ\phi_t\circ \rho)^{-1},$
		for each $t$. 
		\item For any isotopy $\Phi = \{\phi_t\},$ we shall denote by $C(\Phi,\eta)^t$ the smooth function\\ $x\mapsto\eta(\dot \phi_t)(\phi_t(x))$. 
		If  $\Phi = \{\psi_t\}, \Psi = \{\phi_t\}\in Iso_{\eta, \omega}(M)$
		such that $\mathcal{L}_{\dot\phi_t}\eta = \mu_t\eta$, (or $\phi_t^\ast(\eta) = e^{f_t}\eta$),
		for each $t$, then we have $  C(\Phi\circ\Psi, \eta)^t =  C(\Phi, \eta)^t\circ\psi_t + e^{f_t\circ\psi_t}  C(\Psi, \eta)^t,$
		for all $t$. 
		\item So, from the above formula, we derive that if $\Phi = \{\phi_t\},$ is an almost cosymplectic isotopy such that $\phi_t^\ast(\eta) = e^{f_t}\eta$, 
		for each $t$, then, we have  $ C(\Phi^{-1},\eta)^t =  - e^{f_t}  C(\Phi, \eta)^t\circ \phi_t^{-1},$ 
		for all $t$. 	
	\end{enumerate}

\subsection{The $C^0-$topology}\label{SC2}
Let $Homeo(M)$ denote the group of all homeomorphisms of $M$ equipped with the $C^0-$ 
compact-open topology. This is the 
metric topology induced by the following distance
\begin{equation}
d_0(f,h) = \max(d_{C^0}(f,h),d_{C^0}(f^{-1},h^{-1})),
\end{equation}
where $
d_{C^0}(f,h) =\sup_{x\in M}d (h(x),f(x)).$
On the space of all continuous paths $\lambda:[0,1]\rightarrow Homeo(M)$ such that $\lambda(0) = id_M$, 
we consider the $C^0-$topology as the metric topology induced by the following metric  $
\bar{d}(\lambda,\mu) = \max_{t\in [0,1]}d_0(\lambda(t),\mu(t)).$\\

We will need the following facts. 
\begin{remark}\label{Impot-1}
	For any smooth map $\psi$, any vector field $X$, and any differential $k-$form $\beta$, we have $\psi^\ast(\iota(\psi_\ast(X))\beta) = \iota(X)\psi^\ast(\beta).$ Indeed, pick $(k-1)$ vector fields $Y_1,\dots,Y_{k-1}$ on $M$, and compute\\
	$\psi^\ast(\iota(\psi_\ast(X))\beta)(Y_1,\dots,Y_{k-1}) = (\iota(\psi_\ast(X))\beta)(\psi_\ast(Y_1),\dots,\psi_\ast(Y_{k-1})) $
	$$ = \beta\circ\psi (\psi_\ast(X),\psi_\ast(Y_1),\dots,\psi_\ast(Y_{k-1}) ),$$
	$$ = \psi^\ast(\beta)(X, Y_1,\dots,Y_{k-1}),$$
	$$ = (\iota(X)\psi^\ast(\beta))(Y_1,\dots,Y_{k-1}).$$
	In particular, if $\psi$ is a diffeomorphism, then we have 
	$ \iota(\psi_\ast(X))\beta = (\psi^{-1})^\ast(\iota(X)\psi^\ast(\beta)).$
\end{remark}

\begin{remark}\label{Impot-0}
	Let $N_i$, $i = 1, 2$  be two smooth manifolds, and denote by $ \tilde N$ the Cartesian product $N_1\times N_2$. 	
	Let $p_i: \tilde N\rightarrow N_i,$ be the projection for $i = 1, 2$. We call a section of $p_i$, any map $ \sigma_i: N_i\longrightarrow \tilde N$ such that $p_i\circ \sigma_i = id_{N_i}$, $i = 1, 2$.
	\begin{enumerate}
		\item  If  $N_i$ is equipped with a Riemaniann metric $g_i$, $i = 1, 2$, then we equip 
		$\tilde N$ with the product metric $\tilde g$, defined by $\tilde g = p^\ast_1(g_1) + p^\ast_2(g_2)$. 
		\item Let equip $\tilde N$ with the following equivalence relation $\sim$ defined by: $(x, \theta) \sim (y, \beta)$, if and only if, $x = y$. 
		\item Assume that $N_i$ is compact for $i = 1, 2$, then by Tykhonov theorem, so is $\tilde N$. For any smooth function $\tilde F$ on $\tilde N$, let denote by $Cri(\tilde F)$ the set of all its critical points, since $\tilde N $ is compact, then by Morse theory, we have $Cri(\tilde F)\neq \emptyset$, for each smooth function $F$ on $\tilde N$. So, if for each smooth function $F$ on $\tilde N$, we denote by $Cri_\flat(\tilde F)$ the quotient space 
		$Cri(\tilde F)/\sim$, then 
		$$ 1\leq \inf_{F\in C^\infty(\tilde N, \mathbb{R})}\left( \sharp ( Cri_\flat(\tilde F))\right)=: \Gamma (\tilde N),$$ 
		where $\sharp A$ stands for the cardinal of a set $A$. 
		\item It is clear that 
		$$  1\leq \Gamma (\tilde N) \leq \min\{\inf_{F\in C^\infty(N_1, \mathbb{R})}\left( \sharp ( Cri(F))\right), \inf_{H\in C^\infty(N_2, \mathbb{R})}\left( \sharp (Cri(H))\right)\}.$$
		To see the above inequalities, let $H$ be any smooth function on $N_1$, and define a smooth function $F:= H\circ p_1$, on $\tilde N$: Assume that $d_{(a,b)}F = 0$, for some $(a, b)\in \tilde N$, and that $d_aH\neq 0$. Consider the section $S_b: N_1\longrightarrow \tilde N, x\mapsto (x, b)$, of the projection $p_1$, and let 
		$X$ be a vector field on $N_1$ such that $dH(X)(a)\neq 0$ (such a vector field exists since $d_aH\neq 0$ ), and set $Y := (S_b)_\ast(X)$. Since 
		$(a, b)\in \tilde N$ is a critical point for $F$, then 
		$$ 0 = dF(Y)(a,b) = dH_a\left( (d_{S_b(a)}p_1\circ d_aS_b)(X_a) \right) $$ 
		$$  = dH_a\left( (d_a(p_1\circ S_b)(X_a) \right) $$
		$$ = dH_a\left( X_a \right)$$
		$$ = dH(X)(a)\neq 0.$$
		This is a contradiction: Hence, any critical point for $F$, generates a critical point for $H$. This implies that  $ \Gamma (\tilde N)\leq \sharp ( Cri(H))$, for 
		any smooth function $H$ on $N_1$. 
		\item In particular: 
		\begin{itemize}
			\item If $ N_2$ is the unit circle $\mathbb{S}^1$, then $  1\leq \Gamma (\tilde N) \leq 2$, since the hight function on $\mathbb{S}^1$ has two critical points. 
			\item If  $ N_2$ is the interval $[-1, 1]$, then  $   \Gamma (\tilde N)  = 1$, since the function\\ $f:[-1, 1]\longrightarrow \mathbb{R}, t\mapsto t^2$, has 
			one critical point. 
			\item If $ N_2$ is the torus  $\mathbb{T}^{2k}$,  
			then $  1\leq \Gamma (\tilde N) \leq 2k + 1$. 
			
		\end{itemize}

	\end{enumerate}
\end{remark}
\section{Co-Hamiltonian geometry}\label{SC00} 
\begin{proposition}\label{pro4}\cite{T-H-B}
	Let $(M, \eta, \omega)$ be a  cosymplectic manifold. 
	The set $Ham_{\eta, \omega}(M)$ is a Lie group whose Lie algebra is the space $ham_{\eta, \omega}^\ast(M)$.
	
\end{proposition}

\begin{proposition}\label{pro3}\cite{T-H-B}
	Let $(M, \eta, \omega)$ be a  cosymplectic manifold. The set $Ham_{\eta, \omega}(M)$ is a normal subgroup in the group $G_{\eta, \omega}^\ast(M)$. 
	
\end{proposition}

We shall need the following result found in  \cite{MM}.

\begin{lemma}\label{lem-3}(\cite{MM})
	Let $M$ be a manifold and $\eta$, $\omega$ be two differential forms on $M$ with degrees $1$ and $2$ 
	respectively. Consider $\tilde M = M\times \mathbb{S}^1$ equipped with the $2-$form\\ $\tilde\omega: = p^\ast(\omega) + p^\ast(\eta)\wedge \pi^\ast_2(d\theta),$ where $\theta$ is the coordinate
	function on the unit circle $\mathbb{S}^1$,  $p: \tilde M\rightarrow M,$ and $\pi_2: \tilde M\rightarrow \mathbb{S}^1,$ are projection maps. Then,  $(M, \eta, \omega)$ is a  cosymplectic manifold if 
	and only if  $(\tilde M, \tilde\omega)$ be a  symplectic manifold.
\end{lemma}

We have the following facts. 
\begin{corollary}\label{Trasit-0}
Let $(M, \eta, \omega)$ be a cosymplectic manifold, and let $(\tilde M, \tilde\omega)$ be the corresponding  
symplectic manifold. Let $X$ be a vector field on $\tilde M$ such that $ \iota(X)\tilde\omega = dH $, for a certain smooth function on $\tilde M$. 
Then, for any given two sections $S_l:M\ni x\mapsto (x,l)\in\tilde M$, and $S_k:M\ni x\mapsto (x,k)\in\tilde M$ with $l, k\in \mathbb{S}^1$ (fixed), 
of the projection  $p: \tilde M\rightarrow M,$ we have 
$$d(H\circ S_l) -  d(H\circ S_k) =  \left((d\theta)(\iota((\pi_2)_\ast(X))(k) - (d\theta)(\iota((\pi_2)_\ast(X))(l)\right) \eta,$$ 
where $\pi_2: \tilde M\rightarrow \mathbb{S}^1,$ is the projection. In particular, if the manifold $M$ is closed, 
then $ d(H\circ S_l) =  d(H\circ S_k)$. 
\end{corollary} 
{\it Proof.} Since $ \iota(X)\tilde\omega = dH $, we derive that for any $l\in \mathbb{S}^1$ (fixed), we have 
$$ d(H\circ S_l) = S_l^\ast (dH)
$$ 
$$ = S_l^\ast\left( p^\ast(\iota(p_\ast(X))\omega) + p^\ast(\iota(p_\ast(X))\eta)\pi_2^\ast (d\theta) - 
\left( \iota(X)\pi^\ast_2(d\theta) \right) p^\ast(\eta) \right) $$
$$ = \iota(p_\ast(X))\omega - S_l^\ast\left( \iota(X)\pi^\ast_2(d\theta)\right)\eta.$$
So, we derive from Remark \ref{Impot-1}, that $ \iota(X)\pi^\ast_2(d\theta) =  \pi^\ast_2 (\iota((\pi_2)_\ast(X))(d\theta)),$ and hence
$$ d(H\circ S_l)= \iota(p_\ast(X))\omega - S_l^\ast\left(\pi^\ast_2 (\iota((\pi_2)_\ast(X))(d\theta))\right)\eta $$
$$ = \iota(p_\ast(X))\omega - \left((\iota((\pi_2)_\ast(X))(d\theta))\right)(l)\eta.$$ Thus, 
for any given two sections $S_l:M\ni x\mapsto (x,l)\in\tilde M$, and $S_k:M\ni x\mapsto (x,k)\in\tilde M$ with $l, k\in \mathbb{S}^1$ (fixed), 
of the projection map $p: \tilde M\rightarrow M,$ we have 
$$d(H\circ S_l) -  d(H\circ S_k) =  \left((d\theta)(\iota((\pi_2)_\ast(X))(k) - (d\theta)(\iota((\pi_2)_\ast(X))(l)\right) \eta,$$ 
where $\pi_2: \tilde M\rightarrow \mathbb{S}^1,$ is the  projection map.\\ If the manifold $M$ is closed, then we must have 
$\left((d\theta)(\iota((\pi_2)_\ast(X))(k) - (d\theta)(\iota((\pi_2)_\ast(X))(l)\right) = 0$, for all $l, k\in \mathbb{S}^1$: Indeed, if $\left((d\theta)(\iota((\pi_2)_\ast(X))(k) - (d\theta)(\iota((\pi_2)_\ast(X))(l)\right)\neq 0$, for some $l, k$, then 
we shall have $$ \eta = d\left(\frac{1}{Z(k,l)} (H\circ S_l -H\circ S_k) \right) =:df,$$
with $$Z(k,l):=  \frac{1}{\left((d\theta)(\iota((\pi_2)_\ast(X))(k) - (d\theta)(\iota((\pi_2)_\ast(X))(l)\right)}.$$
That is, the $1-$form $\eta$ is exact, and this together with Stokes' theorem would imply that 
\begin{equation}\label{Contr-1}
0 \neq Vol_{\eta, \omega}(M): = \int_M\eta\wedge\omega^n =\int_{ M} d(f\omega^n) =  \int_{\partial M} f\omega^n = 0,
\end{equation}
where $ \dim(M) = 2n +1$.  Relation (\ref{Contr-1}) is a contradiction. $ \blacksquare$
\begin{proposition}\label{Trasit-1}
	Let $(M, \eta, \omega)$ be a compact cosymplectic manifold. 
If $\Phi_F = \{\phi_t\}$ is any co-Hamiltonian (resp. cosymplectic) isotopy such that 
	$I_{\eta, \omega} (\dot\phi_t) = dF_t$, for all $t$, then the isotopy $\tilde\Phi = \{\tilde\phi_t\}$ defined by 
	 $$ \tilde\phi_t: M\times \mathbb{S}^1 \longrightarrow M\times \mathbb{S}^1,$$
	 $$(x, \theta) \mapsto (\phi_t(x), \mathcal{R}_{\Lambda_t(\Phi_F)}(x, \theta) ),$$
	 is a Hamiltonian (resp. symplectic) isotopy of  the symplectic manifold $(\tilde M, \tilde\omega)$,
	 where $$ \mathcal{R}_{\Lambda_t(\Phi)}(x, \theta) = \theta - \int_0^tC(\Phi, \eta)^s\circ\phi_s(x) ds \mod{2\pi}.$$ Conversely, if the map $$ \tilde\phi_t: M\times \mathbb{S}^1 \longrightarrow M\times \mathbb{S}^1,$$
	 $$(x, \theta) \mapsto (\phi_t(x), \mathcal{R}_{\Lambda_t(\Phi_F)}(x, \theta) ),$$ is a Hamiltonian (resp. symplectic) isotopy of  the symplectic manifold $(\tilde M, \tilde\omega)$, then the path $t\mapsto \phi_t$ is a co-Hamiltonian (resp. cosymplectic) isotopy of  $(M, \eta, \omega)$. 
\end{proposition}
{\it Proof.} Assume $\Phi_F = \{\phi_t\}$ to be any cosymplectic isotopy, and consider $p: \tilde M\rightarrow M,$ and $\pi_2: \tilde M\rightarrow \mathbb{S}^1,$ are projection maps. For each $t$, we have $p\circ \tilde{\phi}_t = \phi_t\circ p$
$$\tilde{\phi}_t^\ast(\tilde{\omega}) =  \tilde{\phi}_t^\ast (p^\ast(\omega) + p^\ast(\eta)\wedge \pi^\ast_2(d\theta))$$
$$ =(p\circ\tilde{\phi}_t)^\ast\omega +  (p\circ\tilde{\phi}_t)^\ast\eta\wedge (\pi_2\circ\tilde{\phi}_t)^\ast d\theta$$
$$ =  (\phi_t\circ p)^\ast\omega +  (\phi_t\circ p)^\ast\eta\wedge \pi_2^\ast d\theta$$
$$ =  p^\ast((\phi_t)^\ast\omega) +  p^\ast((\phi_t)^\ast\eta) \wedge \pi_2^\ast d\theta$$
$$ =  p^\ast(\omega) +  p^\ast(\eta) \wedge \pi_2^\ast (d\theta)$$
$$ = \tilde{\omega}.$$
So, $\tilde\Phi = \{\tilde\phi_t\}$ is a symplectic isotopy of  $(\tilde M, \tilde\omega)$. In particular, if 
 $\Phi_F = \{\phi_t\}$ to be any co-Hamiltonian isotopy such that 
 $I_{\eta, \omega} (\dot\phi_t) = dF_t$, for all $t$, then 
 $$ \iota(\dot{\tilde\phi}_t)\tilde\omega = p^\ast(\iota(\dot{\phi}_t)\omega) + p^\ast(\iota(\dot{\phi}_t)\eta)\pi_2^\ast (d\theta) -
 \left( -\eta(\dot\phi_t)\circ p\right) p^\ast(\eta) $$
 $$ = p^\ast\left( \iota(\dot{\phi}_t)\omega  + \eta(\dot\phi_t)\eta\right) +  p^\ast(\iota(\dot{\phi}_t)\eta)\pi_2^\ast (d\theta)$$
 $$ = p^\ast(I_{\eta, \omega} (\dot\phi_t))+  p^\ast(\iota(\dot{\phi}_t)\eta)\pi_2^\ast (d\theta)$$
 $$ = d(F_t\circ p + p^\ast(\iota(\dot{\phi}_t)\eta)\pi_2),$$
 $$ = d\tilde H_t,$$
 with $\tilde H_t:= F_t\circ p + p^\ast(\iota(\dot{\phi}_t)\eta)\pi_2$, i.e., the isotopy $\tilde\Phi = \{\tilde\phi_t\}$ is Hamiltonian, with 
 Hamiltonian $\tilde H_t:= F_t\circ p + p^\ast(\iota(\dot{\phi}_t)\eta)\pi_2$, for each $t$. Conversely, if the map $$ \tilde\phi_t: M\times \mathbb{S}^1 \longrightarrow M\times \mathbb{S}^1,$$
 $$(x, \theta) \mapsto (\phi_t(x), \mathcal{R}_{\Lambda_t(\Phi_F)}(x, \theta) ),$$ is symplectic isotopy of  the symplectic manifold $(\tilde M, \tilde\omega)$, and from $ \tilde{\phi}_t^\ast(\tilde{\omega}) = \tilde \omega$, we derive that 
 $$  p^\ast((\phi_t)^\ast\omega) +  p^\ast((\phi_t)^\ast\eta) \wedge \pi_2^\ast d\theta =  p^\ast(\omega) +  p^\ast(\eta) \wedge \pi_2^\ast (d\theta),$$
 composing the above equality with the pull-back map of any section $S_l: x\mapsto (x,l)$, with $l\in \mathbb{S}^1$ (fixed) of the projection map $p$, we derive that 
 $\phi_t^\ast\omega= \omega$. On the other hand, we derive from $ (\tilde{\phi}_t)_\ast(\frac{\partial}{\partial\theta})  = \frac{\partial}{\partial\theta}$, that 
 $$ -p^\ast(\eta) = \iota(\frac{\partial}{\partial\theta})\tilde\omega.$$ By Remark \ref{Impot-1}, we also have 
$$ \iota(\frac{\partial}{\partial\theta})\tilde\omega = \iota((\tilde{\phi}_t)_\ast(\frac{\partial}{\partial\theta}))\tilde\omega $$ 
$$
 =  (\tilde{\phi}_t^{-1})^\ast ( \iota(\frac{\partial}{\partial\theta})\tilde{\phi}_t^{\ast}(\tilde\omega)) $$ 
 $$
 = (\tilde{\phi}_t^{-1})^\ast ( \iota(\frac{\partial}{\partial\theta})\tilde\omega)$$ 
 $$ = (\tilde{\phi}_t^{-1})^\ast ( - p^\ast(\eta)),$$
 i.e., $ -p^\ast(\eta) =  (\tilde{\phi}_t^{-1})^\ast ( - p^\ast(\eta))$, and  composing this equality with the pull-back map of any section $S_l: x\mapsto (x,l)$, with $l\in \mathbb{S}^1$ (fixed), of the projection map $p$ gives $ -\eta = (\phi_t)^\ast (-\eta)$. Hence, the path $t\mapsto \phi_t$ is a  cosymplectic isotopy of  $(M, \eta, \omega)$. In particular, 
 if $\tilde{\phi}_t$ is a Hamiltonian isotopy of $(\tilde M, \tilde\omega)$ with Hamiltonian $\tilde H$, then from $\iota(\dot{\tilde\phi}_t)\tilde\omega = d\tilde H_t$, we derive that 
 $$ d\tilde H_t = \iota(\dot{\tilde\phi}_t)\tilde\omega $$ 
 $$= p^\ast(I_{\eta, \omega} (\dot\phi_t))+  p^\ast(\iota(\dot{\phi}_t)\eta)\pi_2^\ast (d\theta),$$
 i.e., $p^\ast(I_{\eta, \omega} (\dot\phi_t)) = d\left( \tilde H_t -   p^\ast(\iota(\dot{\phi}_t)\eta)\pi_2\right) $, for each $t$. So, composing the latter equality with any section $S_l: x\mapsto (x,l)$, with $l\in \mathbb{S}^1$, fixed, gives 
 $I_{\eta, \omega} (\dot\phi_t) = d\left( \tilde H_t\circ S_l\right)$, i.e., the path $t\mapsto \phi_t$, is a co-Hamiltonian isotopy 
 with Hamiltonian $ \tilde H\circ S_l$. Note that by Corollary \ref{Trasit-0}, 
 for any given two sections $S_l:M\ni x\mapsto (x,l)\in\tilde M$, and $S_k:M\ni x\mapsto (x,k)\in\tilde M$ with $l, k\in \mathbb{S}^1$ (fixed), 
 of the projection map $p: \tilde M\rightarrow M,$ we have 
 $$d(\tilde H_t\circ S_l) -  d(\tilde H_t\circ S_k) =  \left((d\theta)(\iota((\pi_2)_\ast(\dot{\tilde\phi}_t))(k) - (d\theta)(\iota((\pi_2)_\ast(\dot{\tilde\phi}_t))(l)\right) \eta,$$ 
 $$ = \left(- \eta(\dot{\phi}_t) +  \eta(\dot{\phi}_t) \right)\eta,$$
 $$ = 0,$$
 where $\pi_2: \tilde M\rightarrow \mathbb{S}^1,$ is the projection map. So, the quantity $ d(\tilde H_t\circ S_l) $, uniquely determines 
 a co-Hamiltonian isotopy  $t\mapsto \phi_t$, no matter the choice of the section $S_l:M\ni x\mapsto (x,l)\in\tilde M$, with $l\in \mathbb{S}^1$ (fixed), 
 of the projection  $p: \tilde M\rightarrow M.$ 
  $ \blacksquare$

\subsubsection{Co-Hamiltonian dynamical systems}
\begin{lemma}\label{T-1}
	Let $(M, \eta, \omega)$ be a closed cosymplectic manifold. Then, each co-Hamiltonian diffeomorphism of $(M, \eta, \omega)$  
	must have at least $\Gamma(\tilde M)$  fixed points where $(\tilde M, \tilde \omega)$ is the symplectic manifold $\tilde M = M\times \mathbb{S}^1$ equipped with the symplectic form $\tilde\omega: = p^\ast(\omega) + p^\ast(\eta)\wedge \pi^\ast_2(d\theta),$ where $\theta$ is the coordinate
	function on the unit circle $\mathbb{S}^1$,  $p: \tilde M\rightarrow M,$ and $\pi_2: \tilde M\rightarrow \mathbb{S}^1,$ are projection maps.
\end{lemma}
{\it Proof.}  Let $\psi$ be any co-Hamiltonian diffeomorphism, and let  $\Phi_F = \{\phi_t\}$ be any co-Hamiltonian isotopy such that 
$I_{\eta, \omega} (\dot\phi_t) = dF_t$, for all $t$, with $\phi_1 = \psi$. One defines a Hamiltonian isotopy $\tilde\Phi_F = \{\tilde\phi_t\}$ of  the symplectic manifold $(\tilde M, \tilde\omega)$  as follows: $ \tilde\phi_t: M\times \mathbb{S}^1 \longrightarrow M\times \mathbb{S}^1,$
$(x, \theta) \mapsto (\phi_t(x), \mathcal{R}_{\Lambda_t(\Phi_F)}(x, \theta) ),$
where $ \mathcal{R}_{\Lambda_t(\Phi_F)}(x, \theta) = \theta - \int_0^tC(\Phi_F, \eta)^s(x) ds \mod{2\pi}.$
Furthermore, by Proposition \ref{Trasit-1},  the isotopy $\tilde\Phi_F = \{\tilde\phi_t\}$ is a Hamiltonian 
isotopy of $(\tilde M, \tilde \omega)$. So, by the Arnold conjecture \cite{F-O}, $\tilde\phi_1 = (\psi, \mathcal{R}_{\Lambda_1(\Phi_F)}) $ 
must have at least 
as many fixed points as the minimal number of critical points of a smooth function on the closed symplectic manifold $(\tilde M, \tilde \omega)$. On the other hand, observe that  each fix point $(x_0, l)$ for $\tilde\phi_1 $, generates $x_0$ as a fix point for $\psi$, and nothing is excluding the possibility 
for  $(x_0, s)$ to be another fix point for $\tilde\phi_1 $, for some $s\neq l$. Hence, it is sure from the Arnold conjecture \cite{F-O}, that 
$$\sharp Fix(\psi) \geq \sharp \left( Fix(\tilde\phi_1)/ \sim\right) \geq \Gamma(\tilde M),$$
where $\sim$ is the equivalence relation defined in Remark \ref{Impot-0}. $\blacksquare$ 

\begin{remark}\label{Fix-1}
	The following fact is a consequence of the above lemma. If for any  co-Hamiltonian diffeomorphism $\psi$, we denote by  $Fix(\psi)$ the set of all its critical points, then  the set $ Fix(\psi)$ contains some $z\in M$ such that  $\int_0^1\eta(\dot \phi_s)(z)ds  = 0 $, for any 
	co-Hamiltonian isotopy $\Phi_F = \{\phi_t\}$ with time-one map $ \psi$, and since the function $x\mapsto \eta(\dot \phi_t)(x)$ is constant 
	for each $t$, we then derive that for any 
	co-Hamiltonian isotopy $\Phi_F = \{\phi_t\}$, the function $x\mapsto \int_0^1\eta(\dot \phi_s)\circ\phi_s(x)ds$ is trivial.
\end{remark}

Let $X$ be a co-Hamiltonian vector field of $(M, \eta, \omega)$, and let $\Phi_X$ be its flow. Since 
$I_{\eta, \omega}(X) = d G$, for some $G\in C^\infty(M,\mathbb{R})$, we derive that for each $p\in M$, we have 
$$\frac{d}{ dt}G(\Phi_X^t(z)) = \left( dG(X)\right)(\Phi_X^t(z))$$ 
$$ = \eta(X)^2(\Phi_X^t(z)),$$ $$ \geq 0,$$
for all $t$, and for all $z\in M$. So, not necessary as in symplectic geometry where 
the Hamiltonian is constant along the orbits of its flow; in the cosymplectic case, the first impression we have is that along the 
orbit $t\mapsto\Phi_X^t(p)$, 
the energy function $G$, increases with time, and we have 
\begin{equation}\label{Orbit-1}
G(\Phi_X^t(z)) = G(z)  + \int_0^t\left( \eta(X)^2(\Phi_X^s(z))\right)ds,
\end{equation}
for all $t$, and for all $z\in M$.\\  Therefore, it seems that in such a dynamics system, given a co-Hamiltonian vector field $X$ of $(M, \eta, \omega)$, 
the orbit $t\mapsto\Phi_X^t(p)$ is  $\sigma-$periodic if and only if, $$ \int_0^\sigma\left( \eta(X)^2(\Phi_X^s(p))\right)ds = 0;$$ and since the map 
$s\mapsto \eta(X)^2(\Phi_X^s(p))$, is positive and continuous, then the orbit $t\mapsto\Phi_X^t(p)$ is  $\sigma-$periodic if and only if, 
$ \eta(X)(\Phi_X^s(p)) = 0,$ for all $s\in [0,\sigma]$. 

\begin{proposition}\label{fact}\cite{T-H-B}
	Let $\Phi:= \{\phi_t\}$ be a cosymplectic isotopy, and $\alpha$ be any closed $1-$form of $M$. We have, 
	\begin{equation}\label{fact1}
	\int_M \Delta(\Phi,\eta)\alpha\wedge\omega^n - 	\int_M \Delta(\Phi,\alpha) \eta\wedge\omega^n = n!\langle[\alpha], [\omega^{(n-1)}\wedge \eta]\wedge \widetilde{S}_{\eta, \omega}(\Phi)\rangle,
	\end{equation}
	where $\Delta(\Phi, \bullet) :=   \int_0^1\bullet(\dot{\phi}_t)\circ\phi_t dt $,  $\langle., .\rangle$ is the usual Poincar\'e scalar product, and 
	$\widetilde{S}_{\eta, \omega}$ is the cosymplectic flux homomorphism. 
\end{proposition}
We have the following facts. 

\begin{lemma}\label{prof1} 
	For any  co-Hamiltonian isotopy  $\Phi = \{\phi_t\} $ with time-one map $id_M$ (i.e., a loop at the identity), for each $x\in M$,  the orbit $\mathcal{O}_x^\Phi: t\mapsto \phi_t(x),$ is contractible. 
\end{lemma} 

{\it Proof.} Let $\Phi = \{\phi_t\} $ be any co-Hamiltonian loop. It follows from  Remark \ref{Fix-1} that
 the function $x\mapsto \int_0^1\eta(\dot \phi_t)\circ\phi_t(x)dt$, is trivial. On the other hand, since 
 $\widetilde{S}_{\eta, \omega}(\Phi) = 0$, and the function $x\mapsto\Delta(\Phi,\alpha)(x) $ is constant for each  
 closed $1-$form $\alpha$ of $M$ (see \cite{Tchuiaga2}), we derive from Proposition \ref{fact1} that 
 $$ 0 = \int_M \Delta(\Phi,\alpha) \eta\wedge\omega^n =  \Delta(\Phi,\alpha) \int_M \eta\wedge\omega^n,$$
for all  closed $1-$form $\alpha$ of $M$. Hence,  $ \Delta(\Phi,\alpha) = 0$, for each closed $1-$form $\alpha$ of $M$. $\blacksquare$
\begin{corollary}\label{prof21} 
Let $\Phi = \{\phi_t\} $ be a co-Hamiltonian isotopy, and let   $\alpha$ be a closed $1-$form of $M$. Then,  	
	$$\int_M \Delta(\Phi,\alpha)\eta\wedge\omega^n = 0.$$ In particular, we have 
	$$\min_x\Delta(\Phi,\alpha)(x)\leq 0 \leq \max_x \Delta(\Phi,\alpha)(x).$$ 
\end{corollary}
{\it Proof.} Let $\Phi = \{\phi_t\} $ be any co-Hamiltonian isotopy, and  
define a Hamiltonian isotopy $\tilde\Phi = \{\tilde\phi_t\}$ of  the symplectic manifold $(\tilde M, \tilde\omega)$  as follows: 
$ \tilde\phi_t:= (\phi_t, \mathcal{R}_{\Lambda_t(\Phi)}),$ for each $t$. By the Arnold conjecture \cite{F-O}, 
there exists $(z,s)\in (\tilde M, \tilde\omega)$ such that $ \tilde\phi_1(z,s) = (z, s)$. This implies that, $\int_0^1 C(\Phi, \eta)^t(z)dt = 0$.  Therefore, we derive as in the proof of 
Lemma \ref{prof1} that  $\int_M \Delta(\Phi,\alpha)\eta\wedge\omega^n = 0$, for all  closed $1-$form $\alpha$ of $M$. In particular, from the inequalities, 
$$\min_x\Delta(\Phi,\alpha)(x)\left(\int_M\eta\wedge\omega^n \right) \leq\int_M \Delta(\Phi,\alpha)\eta\wedge\omega^n  \leq \left(\int_M\eta\wedge\omega^n \right)\max_x \Delta(\Phi,\alpha)(x),$$
we derive that $$\min_x\Delta(\Phi,\alpha)(x)\leq 0 \leq \max_x \Delta(\Phi,\alpha)(x),$$
since  $\int_M \Delta(\Phi,\alpha)\eta\wedge\omega^n = 0$, for all  closed $1-$form $\alpha$ of $M$. 
 $\blacksquare$
\begin{lemma}\label{prof22} 
Let $\psi$ be a co-Hamiltonian diffeomorphism. For any $z\in Fix(\psi)$, and for any co-Hamiltonian isotopy
 $\Phi = \{\phi_t\} $ with time-one map $\psi$, the orbit $\mathcal{O}_z^\Phi: t\mapsto \phi_t(z),$ is contractible. 
\end{lemma}
{\it Proof.} Let $\psi$ be a  co-Hamiltonian diffeomorphism, and pick $z\in Fix(\psi)$. Consider $\Phi = \{\phi_t\} $ to be any co-Hamiltonian isotopy with time-one map $\psi$, and  
define a Hamiltonian isotopy $\tilde\Phi = \{\tilde\phi_t\}$ of  the symplectic manifold $(\tilde M, \tilde\omega)$  as follows: 
$ \tilde\phi_t:= (\phi_t, \mathcal{R}_{\Lambda_t(\Phi)}),$ for each $t$. By Remark \ref{Fix-1}, 
the function $x\mapsto \int_0^1\eta(\dot \phi_t)\circ\phi_t(x)dt$ is trivial. This suggests that for each $\theta\in \mathbb{S}^1$, we have 
$$\tilde\phi_1(z,\theta):= (\phi_1(z), \mathcal{R}_{\Lambda_1(\Phi)}(z,\theta)) = (z, \theta),$$
which means that $(z, \theta) $ is a fix point for $\tilde \phi_1$, hence by the contractibility result of the orbits of fixed points for Hamiltonian diffeomorphisms found in \cite{McDuff-SAl}, \cite{Tchuiaga2}, it follows that 
$$0 = \int_{\mathcal{O}^{\tilde\Phi}_{(z,\theta)}}p^\ast(\alpha) = \int_{p(\mathcal{O}^{\tilde\Phi}_{(z,\theta)})}\alpha = 
\int_{\mathcal{O}^{\Phi}_z}\alpha,$$ 
for each $\theta\in \mathbb{S}^1$, and for all  closed $1-$form $\alpha$ of $M$. $\blacksquare$
\begin{lemma}\label{prof27} 
Let $X$ be a co-Hamiltonian vector field such that $I_{\eta, \omega}(X) = d G$. If $\Phi_X$ is the co-Hamiltonian flow generated by $X$,  then 
the function $G$ is constant  along the orbits of $\Phi_X$.  
\end{lemma}
{\it Proof.} For simplicity, set $\Phi_X = \{\phi_t\}$. 
Since $\Phi_X^1$ has a fix point $z\in M$, then we derive  from (\ref{Orbit-1}) that  
$ \int_0^1\left( \eta(X)^2(\phi_s(z))\right)ds = 0$. The map  
$s\mapsto  \eta(X)^2(\phi_s(z)) $ being is positive and continuous, then the vanishing of the integral $ \int_0^1\left( \eta(X)^2(\phi_s(z))\right)ds$,  implies that $ \eta(X)^2(\phi_s(z)) = 0$, for all $s\in [0,1]$. Since the map $x\mapsto \eta(X)^2(\phi_s(x))$ is constant for each $s$, then the latter must be trivial since $M$ is connected, and  $ \eta(X)^2(\phi_s(z)) = 0$, for all $s\in [0,1]$. Therefore, again from (\ref{Orbit-1}), it follows that  $G\circ\Phi_X^s  = G,$
for all $s\in [0,1]$ because the function $x\mapsto \int_0^1\left( \eta(X)^2(\phi_s(x))\right)ds,$ is trivial. $\blacksquare$

\subsection{Co-Hofer geometry}
 The oscillation of any smooth function $f$ is given by the following 
 formula, 
 \begin{equation}
 osc(f) = \max_{x\in M}f(x) - \min_{x\in M}f(x).
 \end{equation}

Observe that any constant function on $M$ admits a trivial oscillation i.e oscillation is degenerated on the space 
$C^\infty(M,\mathbb{R})$. Thus, the oscillation cannot define a norm over the space $C^\infty(M,\mathbb{R})$. 
Its restriction to the space of smooth functions $f$ satisfying $\int_Mf\eta\wedge\omega^n = 0$, is a norm and we have: 
\begin{enumerate}
	\item $0\leq osc(f)$
	\item $osc(f + h)\leq osc(f) + osc(h)$ 
	\item $ osc(\lambda f) = \arrowvert\lambda\arrowvert osc(f)$, $\forall \lambda \in \mathbb{R}$
	\item $osc(f) - osc(h) \leq osc(f + h)$
	\item $osc(f) = 0\Leftrightarrow f = 0$.
\end{enumerate}

\subsubsection{Co-Hofer lengths, \cite{T-H-B}}
Let  $\Phi_F = \{\phi_t\}$ be a co-Hamiltonian isotopy such that $I_{\eta, \omega} (\dot\phi_t) = dF_t$, for all $t$, then in \cite{T-H-B}, its 
two versions of co-Hofer lengths are defined as:
\begin{equation}\label{equ13}
l_{CH}^{(1,\infty)} (\Phi_F) = \int_0^1 \left( osc(F_t) + | C(\Phi_F,\eta)^t|\right) dt,
\end{equation}
and,
\begin{equation}\label{equ14-s}
l_{CH}^{\infty} (\Phi_F) = \max_t \left( osc(F_t) + | C(\Phi_F,\eta)^t|\right). 
\end{equation}
Furthermore, since for any co-Hamiltonian isotopy $\Phi_F = \{\phi_t\}$ we have,\\ $I_{\eta, \omega}(\dot{\Phi}_F^{-t}) = d(- F_t\circ\phi_t)$, for all $t$, and $$\eta(\dot{\Phi}_F^{-t}) = -\eta((\phi_t^{-1})_\ast(\dot{\phi}_t)))  = -\eta(\dot{\phi}_t)\circ\phi_t = -C(\Phi_F, \eta)^t,$$ 
for each $t$,  we derive that 
\begin{equation}\label{equ13-s}
l_{CH}^{(1,\infty)} (\Phi_F^{-1}) = \int_0^1 \left( osc(-F_t) + | -C(\Phi_F,\eta)^t|\right) dt = l_{CH}^{(1,\infty)} (\Phi_F),
\end{equation}
and,
\begin{equation}\label{equ14}
l_{CH}^{\infty} (\Phi_F^{-1}) = \max_t \left( osc(-F_t) + |- C(\Phi_F,\eta)^t|\right) = l_{CH}^{\infty} (\Phi_F). 
\end{equation}
Thus, the lengths $ l_{CH}^{\infty}$, and 
$ l_{CH}^{(1,\infty)}$ are symmetric. Also, for each $\rho\in G_{\eta, \omega}^\ast(M)$, the co-Hamiltonian isotopy
 $ \Psi_\rho:  t\mapsto \rho^{-1}\circ \phi_F^t\circ \rho$ has the same length than $ \Phi_F$.\\ 

Now, consider a symplectic manifold $(\tilde M, \tilde\omega)$, defined as in Lemma \ref{lem-3}.  In \cite{T-H-B}, it is showed
that if $\Phi_F = \{\phi_t\}$ is a co-Hamiltonian isotopy  such that $I_{\eta, \omega} (\dot\phi_t) = dF_t$, for all $t$, with $\phi_1 \neq id_M$, and $\phi_1$ completely displace a closed ball $\mathbf{B}_0 \subset M$, then we have 
\begin{equation}\label{Ene-trans7}
0\textless \frac{1}{4 \pi^2}\int_{\mathbf{C}_0} C_W(\mathbf{B}_0\times \{.\})d\theta\leq  l_{CH}^{(1,\infty)} (\Phi_F),
\end{equation}
for some compact subset $\mathbf{C}_0$ of the unit circle $\mathbb{S}^1$, where  $C_W(\overline{B})$ is the Gromov area (\cite{Lal-McD95}) of a ball $\overline{B}$ on the closed symplectic manifold $(\tilde M, \tilde\omega)$.

\subsubsection{Almost co-Hofer-like lengths, \cite{T-H-B}} 
For any $X\in \chi_{\eta, \omega}(M)$, the closed $1-$form $\iota(X)\omega$ splits as:
\begin{equation}\label{equ5-A}
\iota(X)\omega = \mathcal{H}_\omega + dU_\omega.
\end{equation}
From the above splitting, one defines a norm $\|.\|_{\mathcal{A}C}^\mathcal{S} $ on 
$\chi_{\eta, \omega}(M)$ as follows: For any  $X\in \chi_{\eta, \omega}(M)$, 
\begin{equation}\label{equ7-A}
\|X\|_{\mathcal{A}C}^\mathcal{S} := \|\mathcal{H}_\omega\|_{L^2} + \nu^B(dU_\omega) + \Theta(X) ,
\end{equation}
with
$$ \Theta(X): = \frac{1}{Vol_{\eta, \omega}(M)}|\int_M\eta(X)\eta\wedge\omega^n|,$$
where $ \|.\|_{L^2}$ is the $L^2-$Hodge norm, $\dim(M) = 2n +1$, $Vol_{\eta, \omega}(M): = \int_M\eta\wedge\omega^n $ and $ \nu^B$ is any norm on $ \mathbb{B}^1(M)$ which we assume to be  equivalent to the oscillation norm. \\
Let $ \Phi = \{\phi_t\}\in Iso_{\eta, \omega}(M)$, such that $\mathcal{L}_{\dot\phi_t}\eta = \mu_t\eta$, 
for each $t$, we have  
$$ \|\dot{\phi}_t\|_{\mathcal{A}C} := \|\mathcal{H}_\omega^t\|_{L^2} + osc(U_\omega^t) + \varTheta_t(\Phi),$$
with 
$$ \varTheta_t(\Phi): = \frac{1}{Vol_{\eta, \omega}(M)}|\int_M C(\Phi,\eta)^t\eta\wedge\omega^n|,$$
for each $t$. Therefore, we define the $L^{(1, \infty)}-$version of the almost co-Hofer-like length of $\Phi: =\{\phi_t\}$ as:
\begin{equation}\label{equ8-A}
l_{\mathcal{A}co}^{(1,\infty)} (\Phi) : = \int_0^1\|\dot{\phi}_t\|_{\mathcal{A}C}dt,
\end{equation}
and, $L^{\infty}-$version of the almost co-Hofer-like length of $\Phi$ as:
\begin{equation}\label{equ9-A}
l_{\mathcal{A}co}^{\infty} (\Phi) : = \max_{t\in [0,1]}\|\dot{\phi}_t\|_{\mathcal{A}C}. 
\end{equation}
As in the case of co-Hofer-like length, it seems that in general, we have
\begin{equation}\label{equ10-A}
l_{\mathcal{A}co}^{(1,\infty)} (\Phi) \ne  l_{\mathcal{A}co}^{(1,\infty)} (\Phi^{-1}),
\end{equation}
or,
\begin{equation}\label{equ11-A}
l_{\mathcal{A}co}^{\infty} (\Phi) \ne  l_{\mathcal{A}co}^{\infty} (\Phi^{-1}). 
\end{equation}
\subsubsection{Almost co-Hamiltonian lengths, \cite{T-H-B}}
The restriction of the above lengths to the group $ \mathcal{A}H_{\eta, \omega}(M)$ will be called almost co-Hofer lengths, and denoted $ l_{\mathcal{A}H}^{\infty}$, and 
$ l_{\mathcal{A}H}^{(1,\infty)}$. Indeed, if  $\Phi = \{\phi_t\}$ is an almost co-Hamiltonian isotopy such that 
$ \iota(\dot\phi_t)\omega = dF_t$, and $\mathcal{L}_{\dot{\phi}_t} = \mu_t \eta$ for all $t$, then 
\begin{equation}\label{equ13-A}
l_{\mathcal{A}H}^{(1,\infty)} (\Phi) = \int_0^1 \left( osc(F_t) + \varTheta_t(\Phi)\right) dt,
\end{equation}
and,
\begin{equation}\label{equ14-A}
l_{\mathcal{A}H}^{\infty} (\Phi) = \max_t \left( osc(F_t) + \varTheta_t(\Phi)\right). 
\end{equation}
Note that the lengths $ l_{\mathcal{A}H}^{\infty}$, and 
$ l_{\mathcal{A}H}^{(1,\infty)}$ are symmetric (\cite{T-H-B}).

\section{Approximation and reparameterization of paths}

Given two co-Hamiltonian isotopies $\Phi_F = \{\phi_t\}$ and $\Psi_H = \{\psi_t\}$ such that $I_{\eta, \omega} (\dot\phi_t) = dF_t$, and 
$I_{\eta, \omega} (\dot\psi_t) = dH_t$
for all $t$, one defines the distances between them as: 
\begin{equation}\label{D-equ13}
D_{CH}^{(1,\infty)} (\Phi_F, \Psi_H) = \int_0^1 \left( osc(F_t- H_t) + | C(\Phi_F,\eta)^t -C(\Psi_H,\eta)^t |\right) dt,
\end{equation}
and,
\begin{equation}\label{D-equ14}
D_{CH}^\infty (\Phi_F, \Psi_H) = \max_t\left( osc(F_t- H_t) + | C(\Phi_F,\eta)^t -C(\Psi_H,\eta)^t |\right).
\end{equation}
Similarly, 
given two almost co-Hamiltonian isotopies $\Phi = \{\phi_t\}$ and $\Psi = \{\psi_t\}$ such that $\iota(\dot\phi_t)\omega = dF_t$,  
and 
$ \iota(\dot\psi_t)\omega = dH_t$, 
for all $t$, one defines the distances between them as: 
\begin{equation}\label{I-A}
D_{\mathcal{A}H}^{(1,\infty)} (\Phi, \Psi) = \int_0^1 \left( osc(F_t- H_t) + \frac{1}{Vol_{\eta, \omega}(M)}|\int_M (C(\Phi, \eta)^t -C(\Psi, \eta)^t )\eta\wedge\omega^n|\right) dt,
\end{equation}
and, 
\begin{equation}\label{II-A}
D_{\mathcal{A}H}^\infty (\Phi, \Psi) = \max_t \left( osc(F_t- H_t) + \frac{1}{Vol_{\eta, \omega}(M)}|\int_M (C(\Phi, \eta)^t -C(\Psi, \eta)^t )\eta\wedge\omega^n|\right).
\end{equation}

\subsubsection{Reparameterization of co-Hamiltonian isotopies}\label{RPSC4-1}
We shall need the following basic facts. Let $\Phi_F = \{\phi_t\}$ be any co-Hamiltonian isotopy  such that $I_{\eta, \omega}(\dot\phi_t) = dF_t$, 
for all $t$, and let 
$\zeta :[0,1]\rightarrow[0,1] $ be a smooth function, then the reparameterized path 
$ t\mapsto \phi_{\xi(t)},$ denoted $\Phi^{\zeta},$ is generated by the
element $F^\xi: t\mapsto \dot{\zeta}(t)F_{\xi(t)}$. In fact, set $\psi_t : = \phi_{\zeta(t)}$, and from $\dot{\psi}_t = \dot{\zeta}(t)(\dot{\phi}_{\zeta(t)})$, derive that 
$$ I_{\eta, \omega} (\dot\psi_t) = \dot{\zeta}(t)I_{\eta, \omega}(\dot{\phi}_{\zeta(t)}) = \dot{\zeta}(t)d F_{\xi(t)}  = d(\dot{\zeta}(t)F_{\zeta(t)}),$$
for each $t$. Furthermore, if $\zeta:[0,1]\rightarrow[0,1] $ is a smooth increasing function, then 
\begin{equation}\label{Pa-equ13}
l_{CH}^{(1,\infty)} (\Phi_F^\zeta) = \int_0^1 \left( osc(\dot{\zeta}(t)F_{\zeta(t)}) + |\dot{\zeta}(t)\eta(\dot{\phi}_{\zeta(t)}) |\right) dt,
\end{equation}
$$ = \int_0^1 \left( \dot{\xi}(t)osc(F_{\xi(t)}) + |\eta(\dot{\phi}_{\xi(t)}) |\right) dt$$
$$ = \int_0^1 \left( osc(F_t) + |\eta(\dot{\phi}_t) |\right) dt$$
$$ = l_{CH}^{(1,\infty)} (\Phi_F).$$
Thus, 
\begin{equation}
	l_{CH}^{(1,\infty)} (\Phi_F^\zeta) = l_{CH}^{(1,\infty)} (\Phi_F).
\end{equation}
Also, we have, 
\begin{equation}\label{Pa-equ14}
l_{CH}^{\infty}  (\Phi_F^\zeta) = \max_t\left( osc(\dot{\zeta}(t)F_{\zeta(t)}) + |\dot{\zeta}(t)\eta(\dot{\phi}_{\zeta(t)}) |\right) 
\end{equation}
$$ \leq \left( \max_t(\dot{\zeta}(t))\right)  l_{CH}^{\infty}  (\Phi_F).$$
Thus, 
\begin{equation}
l_{CH}^\infty (\Phi_F^\zeta) \leq  \left( \max_t(\dot{\zeta}(t))\right)l_{CH}^\infty (\Phi_F).
\end{equation}

\subsubsection{Reparameterization of almost co-Hamiltonian isotopies}\label{RPSC4}
Let $\Phi = \{\phi_t\}$ be any almost co-Hamiltonian isotopy   
 such that $\iota(\dot\phi_t)\omega = dF_t$, and \\ $\mathcal{L}_{\dot{\phi}_t}\eta = \mu_t \eta$, 
 for all $t$, and let 
$\zeta :[0,1]\rightarrow[0,1] $ be a smooth function, then the reparameterized path 
$ t\mapsto \phi_{\zeta(t)},$ denoted $\Phi^{\zeta},$ is generated by the
elements $F^\zeta: t\mapsto \dot{\zeta}(t)F_{\xi(t)}$, and\\  $\mathcal{L}_{\dot{\phi}_t^\zeta}\eta = \dot{\zeta}(t)\mu_{\zeta(t)} \eta$. 
Furthermore, if $\zeta:[0,1]\rightarrow[0,1] $ is a smooth increasing function, then 

\begin{equation}
l_{\mathcal{A}H}^{(1,\infty)} (\Phi^\zeta) = l_{CH}^{(1,\infty)} (\Phi).
\end{equation}
Also, we have, 

\begin{equation}
l_{\mathcal{A}H}^\infty (\Phi^\zeta) \leq  \left( \max_t(\dot{\zeta}(t))\right)l_{\mathcal{A}H}^\infty (\Phi).
\end{equation}
\begin{definition}(\cite{Oh-M07}). Given a smooth function $\xi :[0,1]\rightarrow \mathbb{R}$, its norm $\|\xi\|_{ham}$ is defined by $
	\|\xi\|_{ham} = \|\xi\|_{C^0} + \|\dot\xi\|_{L_1},$ 
	with 
	$\|\dot\xi\|_{L_1} = \int_0^1|\dot\xi(t)|dt,$ and \hspace{0.1cm}$\|\xi\|_{C^0} = \sup_t|\xi(t)|.$
\end{definition}
\begin{definition}\label{DefBT-0} 
	Given a smooth function $f\in C^\infty(M\times[0, 1],\mathbb{R} )$, we define its $C^0-$norm 
	as $\lVert f\rVert_{C^0} := \sup_{x, t}|f(x, t) |.$
\end{definition}
\subsection{Boundary flat co-isotopies}\label{SC6} 
\begin{definition}\label{DefBT} If $\Phi_F = \{\phi_t\}$ is any co-Hamiltonian isotopy  such that $I_{\eta, \omega} (\dot\phi_t) = dF_t$, 
	for all $t$, then we shall say that $\Phi_F$ 
	is boundary flat if there exists $\delta\in]0,1[$ such that 
	$F_t = 0,$ for all $t$ in $[0,\delta[\cup]1-\delta,1].$
\end{definition}
\begin{definition}\label{DefBT-A} If $\Phi = \{\phi_t\}$ is any almost co-Hamiltonian isotopy  
	such that $\iota(\dot\phi_t)\omega = dH_t$, and  $\mathcal{L}_{\dot{\phi}_t}\eta = \mu_t \eta$, 
	for all $t$, then we shall say that $\Phi$ 
	is boundary flat if there exists $\delta\in]0,1[$ such that 
	$H_t = 0,$ and $\mu_t = 0 $, for all $t$ in $[0,\delta[\cup]1-\delta,1].$
\end{definition}

For short, a co-Hamiltonian path $ \{ \phi_t\}$ is boundary flat  
if there exists a constant\\ $0\textless\delta\textless1$ such that 
$\phi_t = id_M$ for all $0\leq t\textless\delta,$ and $\phi_t = \phi_1$ for all $1-\delta\textless t\leq 1$, whereas, 
an almost co-Hamiltonian path $ \Phi = \{ \phi_t\}$  such that $\iota(\dot\phi_t)\omega = dH_t$, and  $\mathcal{L}_{\dot{\phi}_t}\eta = \mu_t \eta$, 
for all $t$, is boundary flat  
if there exists a constant $0\textless\delta\textless1$ such that 
\begin{itemize}
\item $H_t = 0,$  and 
\item $\phi_t^\ast(\eta) = \eta$, 
\end{itemize}
for all $t$ in $[0,\delta[\cup]1-\delta,1].$\\ 

We have the following reparameterization lemmas. Their proofs follow from an adaptation of similar results found in \cite{BanTchu}, \cite{Oh-M07}, and \cite{Tchuiaga2}.

\begin{lemma}\label{RL2} If $\Phi_F = \{\phi_t\}$ is any co-Hamiltonian isotopy  such that $I_{\eta, \omega} (\dot\phi_t) = dF_t$, 
	for all $t$, 
	and $\xi_j : [0,1] \rightarrow[0,1],$ $j = 1,2$
	are two smooth monotonic functions, 
	then there exists a positive constant $C(H, \eta)$ which depends on $ H$ and $\eta$, such that,
	$$D_{CH}^{(1,\infty)} (\Phi_F^{\xi_1}, \Phi_F^{\xi_2})\leq C(H, \eta)\|\xi_1 -\xi_2\|_{ham}.$$
\end{lemma}
We shall give a complete proof of Lemma \ref{RL2} later on. 
Here is an immediate consequence of Lemma \ref{RL2}.
\begin{lemma}\label{RL3} Let  $\Phi_{F_i} = \{\phi_t^i\}$ be a Cauchy sequence of co-Hamiltonian isotopies 
	with respect to the metric $D_{CH}^{(1,\infty)} $, 
	and $\xi_l : [0,1] \rightarrow[0,1],$ $l = 1,2$ be two monotonic smooth functions. Let $\epsilon\textgreater0$ be given.
	\begin{enumerate}
\item There exists a positive constant 
$\delta = \delta(F_i)$ which depends on the sequence $\{F_i\}$, and a larger positive integer $j_0 = j_0(F_i)$ which depends on the sequence $\{F_i\}$, 
such that if the inequality 
$\|\xi_1 -\xi_2\|_{ham}\textless\delta$ holds,
then $$D_{CH}^{(1,\infty)} (\Phi_{F_i}^{\xi_1}, \Phi_{F_i}^{\xi_2})\textless\epsilon,$$ 
for all $i\geq j_0$.
\item There exists a positive constant 
$\tau = \tau(F_i)$ which depends on the sequence $\{F_i\}$, and a larger positive integer $j_0 = j_0(F_i)$ which depends on the sequence $\{F_i\}$,
such that  $$D_{CH}^{(1,\infty)} (\Phi_{F_i}\circ\lambda, \Phi_{F_i}\circ \mu)\textless\epsilon,$$ 
for all $i\geq j_0$, whenever $ d_{C^0}(\lambda, \mu) \textless\tau$.
	\end{enumerate} 
	
\end{lemma}
{\it Proof.} For $(1)$, since  $\Phi_{F_i} = \{\phi_t^i\}$ is Cauchy 
with respect to the norm $D_{CH}^{(1,\infty)} $,  
one can choose an integer $j_0$ large enough such that 
$D_{CH}^{(1,\infty)} (\Phi_{F_i}, \Phi_{F_{j_0}})\textless\epsilon/3$ for all $i\geq j_0$. Assume this is done, and 
 apply Lemma \ref{RL2} with $ \Phi_{F_{j_0}}$, $\xi_1$, and $\xi_2$ to derive that 
there exists a constant $C$ which depends on $F_{j_0}$ and $\eta$ such that, 
$
D_{CH}^{(1,\infty)} (\Phi_{F_{j_0}}^{\xi_1}, \Phi_{F_{j_0}}^{\xi_2})
\leq C(F_{j_0}, \eta)\|\xi_1 -\xi_2\|_{ham}.
$
Taking $\delta = \epsilon/3C(F_{j_0}, \eta)$, we compute
$$ D_{CH}^{(1,\infty)} (\Phi_{F_i}^{\xi_1}, \Phi_{F_{i}}^{\xi_2})\leq 
D_{CH}^{(1,\infty)} (\Phi_{F_{i}}^{\xi_1}, \Phi_{F_{j_0}}^{\xi_1}) + D_{CH}^{(1,\infty)} (\Phi_{F_{j_0}}^{\xi_1}, \Phi_{F_{j_0}}^{\xi_2})
+ D_{CH}^{(1,\infty)} (\Phi_{F_{j_0}}^{\xi_2}, \Phi_{F_{i}}^{\xi_2})$$ 
$$\leq \epsilon/3 + \epsilon/3 + \epsilon/3,$$
as long as $\|\xi_1 -\xi_2\|_{ham}\textless\delta$, and $i\geq j_0.$\\ For $(2)$, choose an integer $j_0$ large enough such that 
$D_{CH}^{(1,\infty)} (\Phi_{F_i}, \Phi_{F_{j_0}})\textless\epsilon/3$ for all $i\geq j_0$. From the uniform continuity 
of the maps $ x\mapsto F_{j_0}(x)$, and $x\mapsto C(\Phi_{F_i}, \eta)(x)$ (this last map is a constant map), we derive that there exists $\tau\textgreater 0$ 
such that $\lVert F_{j_0}\circ\lambda -   F_{j_0}\circ\mu\rVert_{C^0} < \frac{\epsilon}{6},$ 
whenever $ d_{C^0}(\lambda, \mu) \textless\tau$. This implies that 
$D_{CH}^{(1,\infty)} (\Phi_{F_{j_0}}\circ\lambda, \Phi_{F_{j_0}}\circ\mu) < \frac{\epsilon}{3},$
whenever $ d_{C^0}(\lambda, \mu) \textless\tau$. Then, we apply the triangle inequality as in $(1)$ to conclude. $\blacksquare$
\begin{lemma}\label{L2} Let $\Phi_F$ be a co-Hamiltonian isotopy, 
	and let $\epsilon$ be a positive real number. Then, 
	there exists a boundary flat co-Hamiltonian isotopy $\Psi_H$  such that 
	\begin{enumerate}
		\item $\Phi_F^0 = \Psi_H^0$, and  $\Phi_F^1 = \Psi_H^1$,
	\item  $D_{CH}^{(1,\infty)} (\Phi_{F}, \Psi_H)\textless\epsilon,$ and 
	\item $\bar{d}(\Psi_H, \Phi_F) \textless \epsilon.$
	\end{enumerate}
\end{lemma}

{\it Proof.} Let $\epsilon$ be a positive real number, and consider $\xi:[0,1]\rightarrow[0,1] $ 
to be any smooth and 
increasing function which is the trivial map $0$ on $[0,  \delta]$ and the constant map $1$ 
on $[1 - \delta, 1]$ where 
$0\textless\delta\textless 1/13$. Therefore, define $ H $ to be the element 
$F^{\xi}$ obtained by re-parametrising  $F$ 
via the curve $\xi$ as explained in Subsection \ref{RPSC4}.   
It follows from the definition of the curve   
$\xi$ that the co-Hamiltonian isotopy $ \Psi_H$ is boundary flat and it 
has the same extremities than $\Phi_F$. 
Applying Lemma \ref{RL2} with $\xi_1 = id_{[0,1]}$ and $\xi_2 = \xi$, we derive 
from the above arguments that 
$ D_{CH}^{(1,\infty)} (\Phi_{F}, \Psi_H) \leq C(F, \eta)\|\xi - id_{[0,1]}\|_{ham},$ 
where $C$ is the constant in Lemma \ref{RL2} which only 
depends on $ F$ and $\eta$. On the other hand, since the maps 
$(t,x)\mapsto\phi_F^{t}(x)$  and $(t,x)\mapsto\phi_F^{-t}(x)$ 
are Lipschitz continuous, then 
there exists a constant $l_0\textgreater0$ which 
depends  only on $F$ such that 
$\bar{d}(\Psi_H, \Phi_F)\leq 
l_0\|\xi - id_{[0,1]}\|_{C^0}\textless l_0\|\xi - id_{[0,1]}\|_{ham}.$ 
Finally, to conclude, it suffices to choose the function $\xi$ so that $\|\xi - id_{[0,1]}\|_{ham}
\leq\min\{\epsilon/C(F, \eta); \epsilon/l_0;\epsilon\}.$ $\blacksquare$\\

{\it Proof of Lemma \ref{RL2}.} Let $\Phi_F$ be a co-Hamiltonian isotopy generated by  $F$.  
\begin{itemize}
	\item Step (1). Consider the normalized function $V = F^{\xi_1} - F^{\xi_2},$ and compute
\end{itemize}
\begin{equation}\label{0refin3}
|V_t| = |\dot\xi_1(t)F_{\xi_1(t)} -\dot\xi_2(t)F_{\xi_2(t)} |
\leq |\dot\xi_1(t)||F_{\xi_1(t)} - F_{\xi_2(t)}| + |\dot\xi_1(t) - \dot\xi_2(t)||F_{\xi_2(t)}|,
\end{equation}
for each $t$. Since the map $(t,x)\mapsto F_t(x)$ is smooth on a compact set, it is Lipschitz, 
i.e., there exists  
a positive constant $k_0$ which depends on $F$ such that 
$\max_{x\in M}|F_{t}(x) - F_{s}(x)| \leq k_0|t - s|$ for all $t,s\in [0,1]$. 
The latter Lipschitz inequality together with (\ref{0refin3}) yields 
\begin{equation}\label{refin3}
0\leq\max_{x\in M}V_t(x)\leq k_0|\dot\xi_1(t)||\xi_1(t) - \xi_2(t)| 
+ \max\{ \max_x (F_{\xi_2(t)}(x)), -\min_x(F_{\xi_2(t)}(x))\}  |\dot\xi_1(t) - \dot\xi_2(t)|.
\end{equation}
Similarly, one can check that 
\begin{equation}\label{refin4}
0\leq-\min_{x\in M}V_t(x)\leq k_0|\dot\xi_1(t)||\xi_1(t) - \xi_2(t)| +
\max\{ \max_x (F_{\xi_2(t)}(x)), -\min_x(F_{\xi_2(t)}(x))\} |\dot\xi_1(t) - \dot\xi_2(t)|.
\end{equation}
Adding (\ref{refin3}) and (\ref{refin4}) member to member, 
and integrating the resulting inequality in the variable $t$ gives 
\begin{equation}\label{a12}
\int_0^1osc(V_t)dt\leq 2k_0\max_t|\xi_1(t) - \xi_2(t)| + 2\max_t(osc(F_{t}))\int_0^1|\dot\xi_1(t) - \dot\xi_2(t)|dt,
\end{equation}
$$ \leq 4\max\{k_0, \max_t(osc(F_{t}))\}\|\xi_1 -\xi_2\|_{ham}.$$
Also, we have $
C(\Phi_F^{\xi_j}, \eta)^t = \eta(\dot\Phi^{\xi_j}(t)) = \dot\xi_j(t)\eta(\dot\phi_{\xi_j(t)}),
$
for each $t$, and hence  
\begin{equation}\label{TS2}
\int_0^1\arrowvert C(\Phi_F^{\xi_1}, \eta)^t - C(\Phi_F^{\xi_2}, \eta)^t\arrowvert dt 
= \int_0^1\arrowvert \dot\xi_1(t)\eta(\dot\phi_{\xi_1(t)}) - \dot\xi_2(t)\eta(\dot\phi_{\xi_2(t)})\arrowvert dt
\end{equation}
$$ \leq \int_0^1\arrowvert \dot\xi_1(t)\eta(\dot\phi_{\xi_1(t)}) - \dot\xi_1(t)\eta(\dot\phi_{\xi_2(t)})\arrowvert dt + \int_0^1\arrowvert \dot\xi_1(t)\eta(\dot\phi_{\xi_2(t)}) - \dot\xi_2(t)\eta(\dot\phi_{\xi_2(t)})\arrowvert dt$$
$$ \leq \int_0^1|\dot\xi_1(t) |\arrowvert \eta(\dot\phi_{\xi_1(t)}) - \eta(\dot\phi_{\xi_2(t)})\arrowvert dt + 
\int_0^1| \eta(\dot\phi_{\xi_2(t)})| \arrowvert \dot\xi_1(t) - \dot\xi_2(t)\arrowvert dt.$$
The Lipschitz nature of the map $t\mapsto \eta(\dot\phi_t)$ implies  
that there exists a positive constant $c_0$ such that 
$|\eta(\dot\phi_t)  - \eta(\dot\phi_s)|\leq c_0|t - s|$ 
for all $s,t\in [0,1]$. 
This Lipschitz inequality together with  (\ref{TS2}) yields	
\begin{equation}\label{refin5}
\int_0^1\arrowvert C(\Phi_F^{\xi_1}, \eta)^t - C(\Phi_F^{\xi_2}, \eta)^t\arrowvert dt
\leq c_0\max_t|\xi_1(t) - \xi_2(t)| + \max_t|C(\Phi_F, \eta)^t |\int_0^1|\dot\xi_1(t) - \dot\xi_2(t)|dt,
\end{equation}
$$ \leq 2\max\{c_0, \max_t|C(\Phi_F, \eta)^t |\}\|\xi_1 -\xi_2\|_{ham}.$$
Finally, adding  (\ref{a12}) and  (\ref{refin5}) member to member implies that
\begin{equation}\label{a12E}
D_{CH}^{(1,\infty)} (\Phi_F^{\xi_1}, \Phi_F^{\xi_2})\leq C(F,\eta)\|\xi_1 -\xi_2\|_{ham},
\end{equation}
where $$C(F, \eta):= 4\max\{\max\{c_0, \max_t|C(\Phi_F, \eta)^t |\}, 2\max\{k_0, \max_t(osc(F_{t}))\} \}.$$ $ \blacksquare$
\begin{lemma}\label{ch2l2}
	Let $\Phi_F = \{\phi_F^t\}$ be a co-Hamiltonian isotopy such that $ \int_M F_t\eta\wedge \omega^n = 0$, for all $t$.  
	Then, there exists a positive constant $C$ which depends only on $F$ and $\eta$ so that, for any given positive real number
	$\epsilon$, one can find a co-Hamiltonian isotopy $ \Phi_H = \{\phi_H^t\}$  such that 
	\begin{enumerate}
		\item $\Phi_H$ is  boundary flat,
		\item $ \int_M H_t\eta\wedge \omega^n = 0$, for all $t$,
		\item $\phi_H^0 = \phi_F^0$ and $\phi_H^1 = \phi_F^1,$
		\item $D_{CH}^\infty (\Phi_F, \Phi_H)\textless 2D_{CH}^\infty (\Phi_F, Id) + C\epsilon,$ 
		 \item in particular, $D_{CH}^\infty (\Phi_H, Id)\textless 3D_{CH}^\infty (\Phi_F, Id) + C\epsilon,$ 
		\item the $C^0-$distance $\bar{d}(\Phi_F, \Phi_H)$ can be made as small as we want,
	\end{enumerate}
	where $Id:t\mapsto id_M$, is the constant path identity. 
\end{lemma}

{\it Proof}. Let $\Phi_F = \{\phi_F^t\}$ be a co-Hamiltonian isotopy, and let $\epsilon$  be a positive real number. 
We choose a smooth reparameterization curve $\chi :[0,1]\rightarrow[0,1]$ with the following properties: 
\begin{enumerate}
	\item $\chi\equiv 0$ near $t = 0$ and $\chi\equiv  1$ near $t = 1$,
	\item $\|\chi - id\|_{C^0}\textless \epsilon,$
	\item $0\leq\chi'(t)\leq 2$, $\forall t\in [0,1]$.
\end{enumerate}
Take $ H := F^{\chi}$: it is easy to see that $ H$ satisfies the three first properties of Lemma \ref{ch2l2}. For $(4)$, by 
the smoothness of $F$ and the compactness of $M$ the map $t\mapsto F_t$ is Lipschitz i.e there exists a constant $K_0$ such that 
$|F_s - F_t|_\infty:=\max_{x\in M}|F_s(x) - F_t(x)| \leq K_0|t - s|,$ for all $s,t\in [0,1].$  
Set $\tilde B_t =  (\chi'(t)F_{\chi(t)} - F_t)$
for all $t$. Since by assumption we have 
$\int_M F_t\eta\wedge \omega^n = 0$, for all $t$, we derive that  $\int_M \tilde B_t\eta\wedge \omega^n = 0,$  $\max_{x\in M}\tilde B_t(x)\geq 0$, and $ -\min_{x\in M}\tilde B_t(x)\geq0,$ for all 
$t$. Compute  
$$ \max_{x\in M}\tilde B_t(x)\leq 2K_0\epsilon + \max\{ \max_x (F_{t}(x)), -\min_x(F_{t}(x))\}(\chi'(t) - 1).$$ Similar consideration applied to $C_t = -\min_{x\in M}\tilde B_t(x)$ implies that
$$-\min_{x\in M}\tilde B_t(x) \leq 2K_0 \epsilon  +\max\{ \max_x (F_{t}(x)), -\min_x(F_{t}(x))\} (\chi'(t) - 1).$$ 
Combining the above statements we obtain :
$$0\leq\max_{x\in M}\tilde B_t (x)-\min_{x\in M}\tilde B_t(x)
\leq 4K_0\epsilon + 2|\chi'(t) - 1|osc(F_t)
\leq 4K_0\epsilon + 2osc(F_t),$$
i.e $$osc(F_t - H_t) = osc(\chi'(t)F_{\chi(t)} - F_t)$$ $$ = osc(\tilde B_t)$$ $$\leq4K_0\epsilon + 2osc(F_t),$$
for all $t$. That is,
\begin{equation}\label{a12-1}
\max_{t}(osc(F_t - H_t))\leq 4K_0\epsilon + 2\max_{t}(osc(F_t)).
\end{equation}
On the other hand, compute,  
\begin{equation}
|C(\Phi_H, \eta)^t - C(\Phi_F, \eta)^t| = |\chi'(t)\eta(\dot{\phi}_{\chi(t)}) -  \eta(\dot{\phi}_{t})|
\end{equation}
$$ \leq
|\chi'(t)\eta(\dot{\phi}_{\chi(t)}) -  \chi'(t)\eta(\dot{\phi}_{t})| + |\chi'(t)\eta(\dot{\phi}_t) -  \eta(\dot{\phi}_{t})|,$$
$$ \leq2L_0|\chi(t) -t| + |\eta(\dot{\phi}_{t}) ||\chi'(t) - 1|,$$
where $L_0$ is the Lipschitz constant of the map $t\mapsto \eta(\dot{\phi}_{t})$.
The above estimates imply that :
\begin{equation}\label{a13}
\max_t|C(\Phi_H, \eta)^t - C(\Phi_F, \eta)^t| \leq 4L_0\epsilon  + \max_t|\eta(\dot{\phi}_{t}) |. 
\end{equation}
We derive from the relations (\ref{a12-1}) and (\ref{a13}) that :
\begin{equation}\label{a16}
D_{CH}^\infty (\Phi_F, \Phi_H)\leq 4(K_0 + L_0)\epsilon + D_{CH}^\infty (\Phi_F, Id).
\end{equation}
For $(5)$,  the maps $(s,x) \mapsto \phi_F^{s}(x)$ and $(s,x) \mapsto (\phi_H^{s})^{-1}(x)$ are Lipschitz continuous 
with respect to metric $\bar d$, we obtain,   
$\bar{d}(\Phi_F, \Phi_H) \leq l_0\|\chi - id\|_{C^0}$ where the constant $l_0\textless\infty$ depends 
only on $F$. $\square$\\

The following results are almost co-Hamiltonian versions of the above lemmas. Their proof can follow (not obviously) by adapting the previous proofs. So, 
shall not prove them. 
\begin{lemma}\label{RL2-A} If $\Phi_H = \{\phi_t\}$  is any almost co-Hamiltonian isotopy such that 
	$\mathcal{L}_{\dot{\phi}_t}\eta = \mu_t \eta$, 
	for all $t$, 
	and $\xi_j : [0,1] \rightarrow[0,1],$ $j = 1,2$
	are two smooth monotonic functions, 
	then there exists a positive constant $(H, \mu)$ which depends on $ H$ and $\mu$, such that,
	$$D_{\mathcal{A}H}^{(1,\infty)} (\Phi_H^{\xi_1}, \Phi_H^{\xi_2})\leq C(H, \mu)\|\xi_1 -\xi_2\|_{ham}.$$
\end{lemma}
 
Here is an immediate consequence of Lemma \ref{RL2-A}.
\begin{lemma}\label{RL3-A} Let  $\Phi_{F_i} = \{\phi_t^i\}$ be a Cauchy sequence of almost co-Hamiltonian isotopies 
	with respect to the metric $D_{\mathcal{A}H}^{(1,\infty)} $, 
	and $\xi_l : [0,1] \rightarrow[0,1],$ $l = 1,2$ be two monotonic smooth functions. Let $\epsilon\textgreater0$ be given.
	\begin{enumerate}
		\item There exists a positive constant 
		$\delta = \delta(F_i)$ which depends on the sequence $\{F_i\}$, and a larger positive integer $j_0 = j_0(F_i)$ which depends on the sequence $\{F_i\}$, 
		such that if the inequality 
		$\|\xi_1 -\xi_2\|_{ham}\textless\delta$ holds,
		then $$D_{\mathcal{A}H}^{(1,\infty)} (\Phi_{F_i}^{\xi_1}, \Phi_{F_i}^{\xi_2})\textless\epsilon,$$ 
		for all $i\geq j_0$.
		\item There exists a positive constant 
		$\tau = \tau(F_i)$ which depends on the sequence $\{F_i\}$, and a larger positive integer $j_0 = j_0(F_i)$ which depends on the sequence $\{F_i\}$,
		such that  $$D_{\mathcal{A}H}^{(1,\infty)} (\Phi_{F_i}\circ\lambda, \Phi_{F_i}\circ \mu)\textless\epsilon,$$ 
		for all $i\geq j_0$, whenever $ d_{C^0}(\lambda, \mu) \textless\tau$.
	\end{enumerate} 
	
\end{lemma}
\begin{lemma}\label{L2-A} Let $\Phi_F$ be an almost co-Hamiltonian isotopy, 
	and let $\epsilon$ be a positive real number. Then, 
	there exists a boundary flat almost co-Hamiltonian isotopy $\Psi_H$  such that 
	\begin{enumerate}
		\item $\Phi_F^0 = \Psi_H^0$, and  $\Phi_F^1 = \Psi_H^1$,
		\item  $D_{\mathcal{A}H}^{(1,\infty)} (\Phi_{F}, \Psi_H)\textless\epsilon,$ and 
		\item $\bar{d}(\Psi_H, \Phi_F) \textless \epsilon.$
	\end{enumerate}
\end{lemma}

\begin{lemma}\label{ch2l2-A}
	Let $\Phi_F = \{\phi_F^t\}$ be an almost co-Hamiltonian isotopy such that $ \int_M F_t\eta\wedge \omega^n = 0$, for all $t$.  
	Then, there exists a positive constant $C$ which depends only on $F$ and $\eta$ so that, for any given positive real number
	$\epsilon$, one can find an almost co-Hamiltonian isotopy $ \Phi_H = \{\phi_H^t\}$  such that 
	\begin{enumerate}
		\item $\Phi_H$ is  boundary flat,
		\item $ \int_M H_t\eta\wedge \omega^n = 0$, for all $t$,
		\item $\phi_H^0 = \phi_F^0$ and $\phi_H^1 = \phi_F^1,$
		\item $D_{\mathcal{A}H}^\infty (\Phi_F, \Phi_H)\textless 2D_{\mathcal{A}H}^\infty (\Phi_F, Id) + C\epsilon,$ 
		\item in particular, $D_{\mathcal{A}H}^\infty (\Phi_H, Id)\textless 3D_{\mathcal{A} H}^\infty (\Phi_F, Id) + C\epsilon,$ 
		\item the $C^0-$distance $\bar{d}(\Phi_F, \Phi_H)$ can be made as small as we want,
	\end{enumerate}
	where $Id:t\mapsto id_M$, is the constant path identity. 
\end{lemma}

 We have the following definitions, each of which characterizes a subset in the group of homeomorphisms.

\begin{definition}
	A  homeomorphism $h$ is called a topological co-Hamiltonian map (or cohameomorphism) in the $L^{(1,\infty)}$ context (resp. in the $L^{\infty}$ context) if,
	\begin{itemize}
		\item there exists $\lambda\in\mathcal{P}(Homeo(M), id)$ with $\lambda(1) = h$, and 
		\item there exists a sequence $\{\Phi_{H_i}\}$ of co-Hamiltonian isotopies which is Cauchy with respect to the metric  
		$D_{CH}^{(1, \infty)}$ 
		(resp. with respect to the metric $D_{CH}^\infty$) such that\\ $\bar d(\Phi_{H_i}, \lambda)\longrightarrow0, i\longrightarrow \infty$.
	\end{itemize}
\end{definition}
We denote by $Cohameo_{\eta, \omega}^{(1,\infty)}(M)$ (resp. $Cohameo_{\eta, \omega}^{\infty}(M)$) the set of all cohameomorphisms  
of any compact cosymplectic manifold $(M,\eta, \omega)$ in the $L^{(1,\infty)}$ context (resp. in the $L^{\infty}$ context).

\begin{definition}
	A  homeomorphism $h$ is called a topological almost co-Hamiltonian map (or almost cohameomorphism) in the $L^{(1,\infty)}$ context (resp. in the $L^{\infty}$ context) if,
	\begin{itemize}
		\item there exists $\lambda\in\mathcal{P}(Homeo(M), id)$ with $\lambda(1) = h$, and 
		\item there exists a sequence $\{\Phi_{H_i}\}$ of almost co-Hamiltonian isotopies which is Cauchy with respect to the metric  
		$D_{\mathcal{A}H}^{(1, \infty)}$ 
		(resp. with respect to the metric $D_{\mathcal{A}H}^\infty$) such that\\ $\bar d(\Phi_{H_i}, \lambda)\longrightarrow0, i\longrightarrow \infty$.
	\end{itemize}
\end{definition}
We denote by $\mathcal{A}cohameo_{\eta, \omega}^{(1,\infty)}(M)$ (resp. $\mathcal{A}cohameo_{\eta, \omega}^{\infty}(M)$) the set of all almost cohameomorphisms  
of any compact cosymplectic manifold $(M,\eta, \omega)$ in the $L^{(1,\infty)}$ context (resp. in the $L^{\infty}$ context).\\

Finally we list the following problems which arise from various definitions introduced in the present paper. These will be subjects of future study.

\begin{enumerate}
	\item In  \cite{Hofer90}, Hofer has proved that for any given two compactly supported Hamiltonian diffeomorphisms $\phi$, and $\psi$ of $\mathbb{R}^{2n}$, on can control the Hofer distance between $\phi$, and $\psi$  by the $C^0-$distance between $\phi$, and $\psi$ whenever $ supp(\phi\circ\psi^{-1})$ is controlled. \\
	Regarding the correlation between cosymplectic manifolds and symplectic manifolds, it could be interesting to investigate the cosymplectic 
	analogue of such a nice result.
\item In \cite{Polt93}, Polterovich showed that one can connect any Hamiltonian diffeomorphism to the identity map via a regular Hamiltonian path.\\ What 
could be the cosymplectic analogue of such an important result?
\item If it happens that any co-Hamiltonian path can be regularized in the sense of  Polterovich \cite{Polt93}, then using the approximations results found in the present paper, one could show that the two norms $\|.\|_{CH}^{(1,\infty)}$, and $\|.\|_{CH}^\infty $ 
are equal.
\item It is know that $Ham_{\eta, \omega}(M)$ is a normal subgroup of $G_{\eta, \omega}^{\ast}(M)$, but regarding some studies performed 
in Banyaga \cite{Ban}, it is tempting to think of the following important 
questions:
\begin{itemize}
\item Is $Ham_{\eta, \omega}(M)$ a simple group?
\item Does $Ham_{\eta, \omega}(M)$  coincide with $[G_{\eta, \omega}^{\ast}(M), G_{\eta, \omega}^{\ast}(M)]$, the commutator subgroup of the group $G_{\eta, \omega}^{\ast}(M)$?
\item Is $Ham_{\eta, \omega}(M)$ a perfect group?
\item Is any cosymplectic  isotopy in $Ham_{\eta, \omega}(M)$ a co-Hamiltonian isotopy?
\end{itemize}
 \item What is the co-Hamiltonian analogues of the Hamiltonian fragmentation properties?
 \item What are the geometric, topological, and algebraic properties of each of the sets 
 $Cohameo_{\eta, \omega}^{(1,\infty)}(M)$, and  $Cohameo_{\eta, \omega}^{\infty}(M)$?
 \item How can one performs canonical measurement in the spaces of cohameomorphisms?
 \item What about the almost cosymplectic versions of the above questions?
 \item 
 The $L^{(1,\infty)}-$version and the $L^{\infty}-$version 
 of the almost co-Hofer energies of $\phi\in \mathcal{A}Ham_{\eta, \omega}(M)$  are respectively are therefore defined by,
 \begin{equation}\label{bny1}
 \|\phi\|_{\mathcal{A}H}^{(1,\infty)} 
 = \inf \left( l_{\mathcal{A}H}^{(1,\infty)}(\Phi) \right), 
 \end{equation}
 and
 \begin{equation}\label{bny2}
 \|\phi\|_{\mathcal{A}H}^\infty 
 = \inf\left(  l_{\mathcal{A}H}^{\infty}(\Phi)\right), 
 \end{equation}
 where each infimum is taken over the set of all almost co-Hamiltonian isotopies 
 $\Phi$ with time-one map equal to $\phi$.  It is clear that the above energies satisfy positivity, triangle inequality, bi-invariance and symmetry. What about the non-degeneracy ?
 The above norms induces two metrics on $ \mathcal{A}Ham_{\eta, \omega}(M)$ as follows: 
\item 
 	Let $ \Phi = \{\phi_i^t\}$ be a sequence of almost co-Hamiltonian isotopies, 
 	let $\Psi = \{\psi^t\}$ be another almost co-Hamiltonian isotopy, and let $\phi : M\rightarrow M$ be a map such that
 	\begin{itemize}
 		\item  $(\phi_i^1)\xrightarrow{C^0}\phi$, and 
 		\item $l_{\mathcal{A}H}^{\infty}(\{\psi^t\}^{-1}\circ\{\phi_i^t\})\rightarrow0,i\rightarrow\infty$.
 	\end{itemize}
 	Then is $\phi = \psi^1$?
 	\item Following \cite{MT}, at any given point $z_0$ of a cosymplectic manifold $(M, \eta, \omega)$ one can attach a cosymplectic orthogonal 
 	complement of a vector subspace $ \Xi \subset T_{z_0}M$ as follows:
 	
 	$$ \Xi^{\perp} := \{Y\in \subset T_{z_0}M: \eta(Y)= 0, \iota(Y)\omega \in \mathcal{I}\Xi\},$$
 	where $ \mathcal{I}\Xi \subset T_{z_0}M^\ast$, is the the annihilator of $\Xi$. 
 	\begin{definition}\cite{MT} Let $(M, \eta, \omega)$ be a cosymplectic manifold. 
 	 A submanifold $ N \subset M$ is called a coisotropic
 		submanifold if:
 		\begin{itemize}
 	\item the Reeb vector field $\xi$ is tangent to $N$, and 
 	\item  $(TN)^{\perp}\subset TN$.
 		\end{itemize}
 	\end{definition}
Now, that we have introduced the notions of cohameomorphisms and almost cohameomorphisms in hand, we can raise the following $C^0-$rigidity questions.
Let $N$ be a coisotropic submanifold of a closed cosymplectic manifold $(M, \eta, \omega)$.  

\begin{itemize}
	\item If $\phi\in Cohameo_{\eta, \omega}^{(1,\infty)}(M)$ is such that $ \phi(N)$ is a smooth submanifold of $M$, then under which condition(s) do we have
	the Reeb vector field $\xi$, tangent to $\phi(N)$, and $(T\phi(N))^{\perp}\subset T\phi(N)$?
		\item If $\psi\in \mathcal{A}cohameo_{\eta, \omega}^{(1,\infty)}(M)$ is such that $ \psi(N)$ is a smooth submanifold of $M$, then under which condition(s) do we have 
		the Reeb vector field $\xi$, tangent to $\psi(N)$, and $(T\psi(N))^{\perp}\subset T\psi(N)$?
		\item If $\phi\in Cohameo_{\eta, \omega}^{(1,\infty)}(M)$ is such that $ \phi(N)$ is a coisotropic submanifold of $M$, then how 
		does $TN$ relate to $T\phi(N)$?
		\item If $\phi\in \mathcal{A}cohameo_{\eta, \omega}^{(1,\infty)}(M)$ is such that $ \psi(N)$ is a coisotropic submanifold of $M$, then how 
		does $TN$ relate to $T\psi(N)$?
\end{itemize}

\end{enumerate}

\end{document}